\documentclass[11pt]{article}
\usepackage{amsmath,amsthm,amssymb,amsfonts}
\usepackage[width=16.5cm,height=23cm]{geometry}
\usepackage{enumerate}
\usepackage{mathrsfs}
\usepackage[cmtip,all]{xy}
\usepackage{tikz}

\usepackage{ifxetex}
\usepackage{enumerate}
\usepackage{mathrsfs}
\usepackage{hyperref}
\usepackage{cleveref}

\theoremstyle{plain}

\newtheorem*{theorem*}{Theorem}

\newtheorem*{proposition*}{Proposition}
\newtheorem*{conjecture*}{Conjecture}

\newtheorem*{corollary*}{Corollary}

\newtheorem{remark*}{Remark}

\newtheorem{defn*}{Definition}

\newcommand{\CC}{\mathbb{C}}

\newcommand{\ZZ}{\mathbb{Z}}
\newcommand{\smsh}{\wedge}

\usepackage{xspace}
\newcommand{\ess}{\mathbb S}

\newcommand{\Gm}{\mathbb G_m}
\newcommand{\Dugger}{the first author\xspace}

\renewcommand{\top}{\mathrm{top}}

\begin{document}

\title{The Multiplicative Structures on Motivic Homotopy Groups}
\author{Daniel Dugger 
\and 
Bj{\o}rn Ian Dundas 
\thanks{The second author was supported by RCN project no. 312472 “Equations in Motivic Homotopy." 
He wishes to thank the Department of Mathematics F. Enriques, 
University of Milan,
for hospitality and support.}
\and
Daniel C.\ Isaksen
\thanks{The third author was supported by National Science Foundation grant DMS-2202267.  
He wishes to thank Eva Belmont for consultations on 3-primary stable homotopy groups.} 
\and
Paul Arne {\O}stv{\ae}r    
\thanks{The fourth author was supported by RCN project no. 312472 “Equations in Motivic Homotopy"
and a Guest Professorship awarded by The Radboud Excellence Initiative.} 
}

\date{December 2, 2022}

\maketitle

\begin{abstract}We reconcile the multiplications on the homotopy rings of motivic ring spectra used by Voevodsky and Dugger.  {\color{black} While the connection is elementary and similar phenomena have been observed in situations like supersymmetry, neither we nor other researchers we consulted were aware of the conflicting definitions and the potential consequences.  Hence this short note.}
\end{abstract}
\noindent{Mathematical Subject Classification (2010): 14F42 (primary), 13A02 (secondary).}
\vspace{\baselineskip}

The homotopy groups of a motivic spectrum $E$ form a $\ZZ\times\ZZ$-graded abelian group $\pi_{*,\star}E$.   If $E$ is a motivic ring spectrum, then the multiplication induces a ring structure on $\pi_{*,\star}E$, which, if $E$ is commutative, should be graded commutative, as explained by \Dugger in \cite{MR3180827}.  In \cite{MR2031198} Voevodsky displays the dual Steenrod algebra $\mathcal A_{*,\star}$ as a ring with graded commutativity $x\cdot y=y\cdot x\cdot(-1)^{ac}$ for $x\in\mathcal A_{a,b}$ and $y\in\mathcal A_{c,d}$ 
-- the same convention is used in \cite{MR3730515} and \cite{zbMATH07015021} -- 
while 
\cite{MR3180827} yields $x\cdot y=y\cdot x\cdot(-1)^{(a-b)(c-d)}\cdot(-1)^{bd}$.  These are different formulas: for instance, Voevodsky claims $\tau_0\tau=\tau\tau_0$ and according to \cite{MR3180827} we must have that $\tau_0\tau=-\tau\tau_0$.

The authors were distressed to discover this, and worryingly enough, none of those we consulted had discovered the discrepancy (although \cite{MR3180827} claims that the Betti realization is not a ring map).  Was there a subtle mistake buried in the literature somewhere?
Something was surely wrong.  But what?


\section*{Don't panic}
\label{sec:dontpanic}
Fortunately, the results are not irreconcilable, and in fact the solution is already to be found in \cite[Proposition 7.2]{MR3180827}:
\begin{quote}
  ``the'' homotopy ring of a motivic ring spectrum $A$ is not canonical.  
\end{quote}
Let us recall the outline of this story:

\begin{enumerate}[(1)]
\item
Taking as given the usual bigraded family of spheres $S^{p,q}$, one obtains a bigraded abelian group $\pi_{*,\star}A=\bigoplus_{p,q} \pi_{p,q}A$.  But equipping this with a product requires fixing a choice of isomorphisms $\phi_{a,b}\colon S^{a_1,a_2}\wedge S^{b_1,b_2}\cong S^{a_1+b_1,a_2+b_2}$ in the stable homotopy category.
For the product to be associative, a set of familiar pentagonal diagrams has to commute; when this happens let us say that the collection of $\phi$-isomorphisms is {\bf coherent}.
\item Let $\ess$ denote the motivic sphere spectrum. 
The set of  coherent collections of $\phi$-isomorphisms
is a torsor for the group $Z^2(\ZZ\times \ZZ,(\pi_{0,0} \ess)^*)$ of reduced $2$-cocycles on the group $\ZZ\times \ZZ$ with values in the group of units in the ring $\pi_{0,0}\ess$.  In other words, if we fix one collection of coherent $\phi$-isomorphisms then any other such collection differs from it by such a reduced 2-cocycle.   Recall here that a function $\alpha\colon \ZZ\times \ZZ\rightarrow (\pi_{0,0}\ess)^*$ is a 2-cocycle when $\alpha(u+v,w)\cdot \alpha(u,v)=\alpha(v,w)\cdot \alpha(u,v+w)$ for $u,v,w\in \ZZ^2$, and is reduced when $\alpha(0,0)=1$.  

\item Two different choices of coherent $\phi$-isomorphisms typically lead to two different ring structures on $\pi_{*,\star}A$.   The difference 2-cocycle is a coboundary precisely when there is a bigraded isomorphism between these rings that multiplies elements of each bidegree $a=(a_1,a_2)$ by a fixed unit $e_{a}\in \pi_{0,0}(\ess)^*$.  Such isomorphisms are called {\bf standard isomorphisms} in \cite{MR3180827}.
\end{enumerate}

See \cite[Section 7]{MR3180827} for details on the above.

It turns out that the $\phi$-isomorphisms chosen by the first author in \cite{MR3180827} lead to a different ring structure on $\pi_{*,\star}A$ than the one used by Voevodsky, even up to standard isomorphism.  Of course we can still translate between the two rings, and it is not exactly that one choice is right and one is wrong---if a person keeps their wits about them as far as remembering the different conventions, there are no contradictions.
But below we will analyze a collection of different choices and make some suggestions about which ones seem ideal.  We stress that the underlying symmetric monoidal structure of motivic spectra and the definition of homotopy groups are the same in \cite{MR3180827} and \cite{MR2031198}, it is only the choice of coherent $\phi$-isomorphisms (not explicitly spelled out in \cite{MR2031198}, but in some sense there implicitly) that differs.

That multigraded objects have flexibility in sign conventions has been observed in situations other than motivic homotopy theory, for instance, in supersymmetry  \cite{MR1701597}.  We comment on this, as well as on the connection to equivariant theory, in Remarks~\ref{remark:equivariant} and ~\ref{remark:supersymmetry} below.

\section*{The signs they are a-changin'}
 
 Regardless of the base scheme, $\pi_{0,0}\ess$ always contains the following four (not necessarily distinct) square roots of $1$: $1,-1,\epsilon$ and $-\epsilon$, where $-1$ and $\epsilon$ are given by $g\mapsto g^{-1}$ on the topological and Tate circles, $S^1$ and $\Gm$, respectively.
  When choosing the coherent isomorphisms
  $$S^{a_1,a_2}\wedge S^{b_1,b_2}\cong S^{a_1+b_1,a_2+b_2},$$
  where $S^{a_1,a_2}=(S^1)^{\smsh (a_1-a_2)}\smsh\Gm^{\smsh a_2}$, the convention in \cite{MR3180827} was as follows:
every time two $S^1$'s are moved past each other, the sign $-1$ appears, and every time two $\Gm$'s are moved past each other, we get an $\epsilon$.  But swapping $S^1$'s and $\Gm$'s is \emph{not} assigned any punishment in \cite{MR3180827}.  This convention makes sense if $S^1$ and $\Gm$ are regarded as generic objects without any special relation between them, which was the case in the more general settings treated in \cite{MR3180827}.  However, this particular choice raises a problem: when the ground field is the complex numbers  Betti realization sends $\Gm$ to $\Gm(\CC)\simeq S^1$, so moving a $\Gm$ past an $S^1$  \emph{is} detected in topology.  Consequently, with these conventions the Betti realization map $\pi_{*,\star}X\rightarrow \pi_*X(\CC)$ is not a ring homomorphism---there is an annoying sign that comes up (cf. \cite[Proposition 1.19]{MR3180827}).

A better approach is to recognize that the isomorphism $S^{a_1,a_2}\smsh S^{b_1,b_2}\cong S^{a_1+b_1,a_2+b_2}$ should involve $a_2(b_1-b_2)$ swaps of $\Gm$'s past $S^1$'s and we can choose to include a ``generalized sign'' factor to track this.  To this end, choose once and for all a unit $u\in\pi_{0,0}\ess$.  In our applications we will have $u^2=1$ and $u$ will play the role of a ``generalized sign'', but the basic setup only needs $u$ to be invertible.  If $A$ is any motivic ring spectrum with unit map $\eta\colon \ess\to A$ we may consider the $\ZZ\times\ZZ$ -graded ring  $(\pi_{*,\star}A,\cdot)$ provided by \cite{MR3180827} and we may consider the alternative $(\pi_{*,\star}A,\cdot_u)$ with $$x\cdot_u y=x\cdot y\cdot \eta u^{a_2(b_1-b_2)}$$ when $x\in\pi_{a_1,a_2}A$ and $y\in\pi_{b_1,b_2}A$ (``punishing'' each swap of $\Gm$'s past $S^1$'s by multiplying with $u$).
Here $\alpha_u((a_1,a_2),(b_1,b_2))=\eta u^{a_2(b_1-b_2)}$ is the 2-cocycle from our story.  
The cocycle condition gives associativity of $\cdot_u$, and the other axioms for a ring follow readily.  If $A$ is commutative then the same proof as for \cite[Proposition 1.18]{MR3180827} shows that  $x\cdot y=y\cdot x\cdot (-1)^{(a_1-a_2)(b_1-b_2)}\epsilon^{a_2b_2}$. So
\begin{align*}
  x\cdot_u y&=y\cdot_u x\cdot (-1)^{(a_1-a_2)(b_1-b_2)}\eta(\epsilon^{a_2b_2}u^{a_2(b_1-b_2)}u^{-b_2(a_1-a_2)})\\
  &=y\cdot_u x\cdot (-1)^{(a_1b_1+a_1b_2+a_2b_1+a_2b_2)}\eta
  (\epsilon^{a_2b_2}u^{a_2b_1-a_1b_2})\\
  &=y\cdot_u x\cdot (-1)^{a_1b_1}\eta(-u)^{a_2b_1-a_1b_2}\eta(-\epsilon)^{a_2b_2}.
\end{align*}
In particular, when $\eta(\epsilon)=\eta(u)=-1$ then $\eta(-\epsilon)=\eta(-u)=1$ and thus
$$x\cdot_u y=y\cdot_u x\cdot (-1)^{a_1b_1}.$$
This is exactly Voevodsky's convention for commutativity in the dual Steenrod algebra:  graded commutativity with respect to the total grading (cf. \cite[Theorem 2.2]{MR2031198}).

\begin{remark*}
We used a special $2$-cocycle in the above computations, but this wasn't necessary.  
For any reduced 2-cocycle $\alpha$ we can define $x\cdot_\alpha y=x\cdot y\cdot \alpha((a_1,a_2),(b_1,b_2))$, and then there is an associated commutativity formula of the form
\[ x\cdot_u y=y\cdot_u x \cdot w\left ((a_1,a_2),(b_1,b_2) \right) \]
where $w$ is a 2-cocycle that is skew-symmetric in the sense of $w(a,b)=w(b,a)^{-1}$.  In fact
\[ w(a,b)=(-1)^{(a_1-a_2)(b_1-b_2)}\epsilon^{a_2b_2} \alpha(a,b)^{-1}\alpha(b,a).
\]
\end{remark*}

\begin{proposition*}
  An invertible element $u$ in $\pi_{0,0}\ess$ gives a functor $A\mapsto(\pi_{*,\star}A,\cdot_u)$ from motivic ring spectra to $\ZZ\times\ZZ$-graded rings.

  Choosing $u=-1$ or $u=\epsilon$ gives graded rings conforming with Voevodsky's commutativity formulas for ring spectra $A$ having the property that $\eta_A(\epsilon)=-1$. Choosing $u=1$ gives the multiplication in \cite{MR3180827}. 

  Also, over the complex numbers, when choosing $u=-1$ or $u=\epsilon$ Betti realization gives a map of (commutative) graded rings by forgetting weight.
\end{proposition*}

For a given choice of $u$, we can ask whether the rings $(\pi_{*,\star}A,\cdot)$ and $(\pi_{*,\star}A,\cdot u)$ happen to be isomorphic via a standard isomorphism.  Deciding this is equivalent to checking whether $\alpha_u$ is a coboundary.  But if $\beta$ is a $1$-cochain then $(\delta\beta)(a,b)=\beta(a)-\beta(a+b)+\beta(b)$ and is therefore  symmetric in $a$ and $b$.  As $\alpha_u$ is not symmetric, it is not a coboundary.  

If $u^2=1$ then the subgroup $\langle u\rangle$ of $(\pi_{0,0}\ess)^*$ is just $\ZZ/2$, and since $B(\ZZ\times \ZZ)$ is the $2$-torus we have $H^2(\ZZ\times \ZZ;\ZZ/2)=\ZZ/2$.  So as far as twisting by $u$ goes (once $u$ is fixed), there are only two different standard isomorphism classes of homotopy rings that can arise: these are represented by the products $\cdot$ and $\cdot_u$ that we saw above.
 Allowing arbitrary twists from the subgroup $\{1,-1,\epsilon,-\epsilon\}$ increases the number of possibilities to four.

\begin{remark*}
    \label{remark:equivariant}
    Of course,
    these considerations hold in situations other than motivic homotopy theory.  An interesting example is that of $C_2$-equivariant spectra.  When over the real numbers, evaluating at complex points gives a symmetric monoidal functor from motivic spectra to $C_2$-equivariant spectra, where $S^1(\CC)$ corresponds to the trivial representation and $\Gm(\CC)$ to the sign representation $\sigma$.  Thus, choosing your $u$'s in the same way in the motivic and in the $C_2$-equivariant setting gives that Betti realization induces a map of (commutative) bigraded rings.

In this $C_2$-equivariant context, in addition to the forgetful map to nonequivariant spectra there is also the fixed-point functor $A\mapsto \phi A$.  This induces maps of groups $\pi_{p,q}(A)\rightarrow \pi_{p-q}(\phi A)$, and so we can ask whether $\pi_{*,\star}(A)\rightarrow \pi_*(\phi A)$ is a ring homomorphism.  For $u=-1$ it is not, but for $u=\epsilon$ it is.  For this reason we suggest that $u=\epsilon$ is the best choice for both motivic and $C_2$-equivariant homotopy.  With this convention the graded-commutativity formula for the homotopy ring of a ring spectrum is
\[ xy=yx\cdot (-1)^{a_1b_1}(-\epsilon)^{a_2b_1+a_1b_2+a_2b_2}=yx\cdot (-\epsilon)^{(a_1-a_2)(b_1-b_2)}\cdot \epsilon^{a_1b_1},\]
for $x\in \pi_{a_1,a_2}A$ and $y\in \pi_{b_1,b_2}A$.
\end{remark*}

\begin{remark*}
\label{remark:supersymmetry}
We end by  mentioning the connections to supersymmetry.  
  Choosing $u=1$ corresponds to the ``Deligne convention'' (see \cite[1.2.8]{MR1701597}) where commuting something in degree $a+b\sigma$ with something in degree $c+d\sigma$ would introduce the penalty $(-1)^{ac+bd\sigma}$ (where $(-1)^\sigma$ is the twist on the sign representation)
  while choosing $u=\epsilon$ would result in the ``Bernstein convention'' with sign
  $(-1)^{(a+b\sigma)(c+d\sigma)}$.
\end{remark*}

\section*{The effect on the motivic stable homotopy ring}

These considerations led us to wonder whether any well-known
relations in the motivic stable homotopy groups
change under the different sign conventions.
For example, do any of the relations in \cite{MR3141814}
depend on the sign convention?
See \cite{zbMATH07303324} for a recent survey article on motivic stable homotopy groups.

First of all, the relation $(1-\epsilon)\eta^2 = 0$
witnesses that $\epsilon$ \emph{must} play a role
in graded commutativity.  When we commute the element
$\eta$ in $\pi_{1,1}$ past itself, a factor of
$\epsilon$ appears.  Note that $2 \eta^2$ is not zero
in general; this is detected in the 
$\mathbb{R}$-motivic homotopy groups.

Consider the list
\[
\rho, \eta, \nu, \sigma, \eta_\top, \nu_\top, \sigma_\top
\]
of elements of degrees 
\[
(-1,-1), (1,1), (3,2), (7,4), (1,0), (3,0), (7,0)
\]
respectively.  These seven elements are defined in the
motivic stable homotopy ring over any base.  As far as 
we are aware, the only way to produce additional ``universal"
examples is to assemble these elements with Toda brackets.

By inspection, it turns out the commutativity
relations amongst these elements are the same
when $u = 1$ or $u = \epsilon$.
The ``error" factor $\epsilon^{a_2 b_1 + a_1 b_2}$ is not
equal to one in some cases.  However, in all such cases,
we are saved by the relations $(1-\epsilon)\rho = 0$
and $(1-\epsilon)\eta = 0$.

This observation led us to search further for an explicit example where the
cases $u = 1$ and $u =\epsilon$ give different
commutativity relations in the motivic stable
homotopy ring.
We inspected the $2$-complete $\mathbb{R}$-motivic 
stable homotopy ring in a large range \cite{MR4461846},
and we found no possible differences.  
Similarly, a brief, speculative investigation of
$3$-complete homotopy yielded no examples.

On the other hand, assume that $\tau$ detects a stable homotopy
element of degree $(0,-1)$.  This assumption holds, for example,
in the $p$-complete context over the field $\mathbb{C}$.
Then the cases $u = 1$ and $u = \epsilon$ give different commutativity
relations.  For example, if $u = 1$, then $\tau \nu = - \nu \tau$
in $\pi_{3,1}$; but if $u = \epsilon$, then $\tau \nu = - \epsilon \nu \tau$.

This investigation
led us to notice a pattern in the $2$-complete
$\mathbb{R}$-motivic stable homotopy groups that had not been
previously observed.

\begin{conjecture*}
Let $\alpha$ have degree $(s,w)$ in the $2$-complete
$\mathbb{R}$-motivic stable homotopy ring.  If $(1-\epsilon) \alpha$
is non-zero, then $w$ is even.
\end{conjecture*}

\begin{footnotesize}
\bibliographystyle{plain}
\bibliography{mhh}
\end{footnotesize}
\vspace{0.1in}

\begin{center}
Department of Mathematics, University of Oregon, USA\\
email: ddugger@uoregon.edu
\end{center}

\begin{center}
Department of Mathematics, University of Bergen, Norway\\
email: dundas@math.uib.no
\end{center}
\begin{center}
Department of Mathematics, Wayne State University, Detroit, Michigan, USA\\
email: isaksen@wayne.edu
\end{center}

\begin{center}
Department of Mathematics F. Enriques, University of Milan, Italy \\ 
Department of Mathematics, University of Oslo, Norway \\
email: paul.oestvaer@unimi.it, paularne@math.uio.no
\end{center}

\end{document}